# A Lemma and a Conjecture on the Cost of Rearrangements

*Alberto Bressan*

S.I.S.S.A. - Via Beirut 4, Trieste 34014, Italy
e-mail: bressan@sissa.it

## 1 - The main result

Consider a stack of books, containing both white and black books. Suppose that we want to sort them out, putting the white books on the right, and the black books on the left (fig. 1). This will be done by a finite sequence of elementary transpositions. In other words, if we have a stack of all black books of length $a$ followed by a stack of all white books of length $b$, we are allowed to reverse their order at the cost of $a + b$. We are interested in a lower bound on the total cost of the rearrangement. Assume that initially the white and black books are highly mixed. By this we mean that inside every segment of length $\varepsilon$ one can find at least $\kappa\varepsilon$ white books and at least $\kappa\varepsilon$ black books, for a given $\kappa \in\, ]0,1[$. We want to show that, as $\varepsilon \to 0$, the cost of sorting the books grows like $|\log \varepsilon|$.

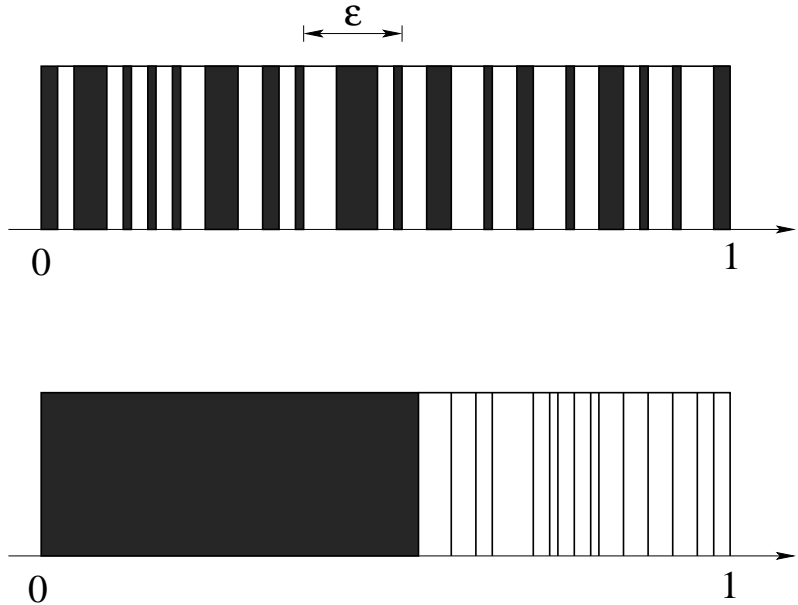

figure 1

We reformulate the problem in more precise terms. Consider the set $\mathcal{F}$ of all piecewise constant functions $f : [0,1] \mapsto \{0, 1\}$. We think of $\mathcal{B} \doteq \{x\,;\ f(x) = 1\}$ as the set of black books. We say that $g$ is obtained from $f$ by an *elementary transposition* if the following holds. For some $y \in [0, 1[$ and $a, b > 0$, one has

$$f(x) = \begin{cases} 1 & \text{if } y < x < y + a, \\ 0 & \text{if } y + a < x < y + a + b, \end{cases}$$

$$g(x) = \begin{cases} 0 & \text{if } y < x < y+b, \\ 1 & \text{if } y+b < x < y+b+a, \\ f(x) & \text{if } x \notin [y,\, y+a+b]. \end{cases}$$

In this case, we define the *cost* of the transposition as $a+b$. By a *rearrangement* of the function $f$ we mean a sequence $f_0, f_1, \ldots, f_n \in \mathcal{F}$ such that

(i) $f_0 = f$,

(ii) $f_n$ is nondecreasing,

(iii) each $f_j$ is obtained from $f_{j-1}$ by an elementary transposition.

The function $f_n$ is then characterized as

$$f_n(x) = \begin{cases} 0 & \text{if } 0 < x < 1 - \xi, \\ 1 & \text{if } 1 - \xi < x < 1, \end{cases} \qquad \xi \doteq \int_0^1 f(s)\,ds. \tag{1.1}$$

The *cost of the rearrangement* $\Gamma(f_0, f_1, \ldots, f_n)$ is of course defined as the sum of the costs of the transpositions.

Now let $0 < \kappa < 1$ be given. We say that the initial configuration $f$ is *well stirred up to scale* $\varepsilon$ provided that

$$\kappa\varepsilon < \int_y^{y+\varepsilon} f(x)\,dx \leq (1-\kappa)\varepsilon \qquad \text{for all } y \in [0,\, 1-\varepsilon]. \tag{1.2}$$

We now show that, as $\varepsilon \to 0$, the cost of a rearrangement grows with logarithmic rate.

**Lemma.** *Assume that the initial configuration $f$ is well stirred up to scale $\varepsilon$. Then the cost of any rearrangement satisfies*
$$\Gamma(f_0, \ldots, f_n) \geq C_\kappa |\log \varepsilon| \tag{1.3}$$
*for some constant $C_\kappa$ independent of $f, \varepsilon$.*

**Proof.** For each $s \in [0, \kappa]$, define $V(s)$ as the minimum cost of any sequence of transpositions $f_0, f_1, \ldots, f_\nu$ such that $f_0 = f$ and $f_\nu \equiv 1$ on some interval with length $\geq s$. By (1.1)-(1.2) it follows

$$\Gamma(f_0, \ldots, f_n) \geq V(1-\xi) \geq V(\kappa/2), \tag{1.4}$$

$$V(s) \geq s - \varepsilon, \tag{1.5}$$

$$V(s) \geq \min_{0 < \sigma < s} V(s-\sigma) + V(\sigma) + \kappa^2 s + (1-\kappa^2)\sigma \qquad s > \varepsilon. \tag{1.6}$$

The last inequality is due to the fact that, in order to create a stack of white books of length $\geq s$, we need first to form a stack of white books with length $s - \sigma$, removing toward the right all the black books that were initially located within them. The amount of such black books can be estimated as

$$\geq \kappa\big((s-\sigma) - (1-\kappa)\varepsilon\big) \geq \kappa^2(s-\sigma)$$

for $s - \sigma \geq \varepsilon$. Then we need to form a second stack of white books of length $\sigma$, and finally join together the two stacks, switching the second white stack with the black stack in the middle, at the cost of $\sigma + \kappa^2(s-\sigma)$. Notice that the possible choice $s - \sigma < \varepsilon$, not considered in the previous argument, certainly does not attain the minimum. Indeed, its cost is

$$\geq V(s-\varepsilon) + (1-\kappa^2)(s-\varepsilon).$$

Using (1.5)-(1.6) one finds
$$V(s) \geq (1+\kappa^2)s - 2\varepsilon.$$

Continuing by induction, for every integer $n$ we obtain

$$V(s) \geq (1+n\kappa^2)s - 2^n\varepsilon.$$

For each $\varepsilon > 0$ we now define $n(\varepsilon)$ as the largest integer $n$ such that

$$(1+n\kappa^2)\frac{\kappa}{2} \geq 2^{n+1}\varepsilon. \tag{1.7}$$

Of course, this implies

$$\Gamma(f_0, \ldots, f_n) \geq V(\kappa/2) \geq \frac{\kappa^3}{4} n(\varepsilon).$$

By (1.7), an elementary estimate yields $n(\varepsilon) \geq C'_\kappa |\log \varepsilon|$, for some constant $C'_\kappa$. This proves the lemma.

## 2 - A conjecture on mixing vector fields

It is interesting to speculate whether a continuous version of the above lemma can be true, for mixing vector fields on a compact manifold. To fix the ideas, consider the two-dimensional torus $K \doteq \mathbb{R}^2/\mathbb{Z}^2$, with coordinates $x = (x_1, x_2) \in [0,1[ \times [0,1[$. Consider the set (fig. 2)

$$A \doteq \{(x_1, x_2); \ 0 \leq x_2 \leq 1/2\} \subset K.$$

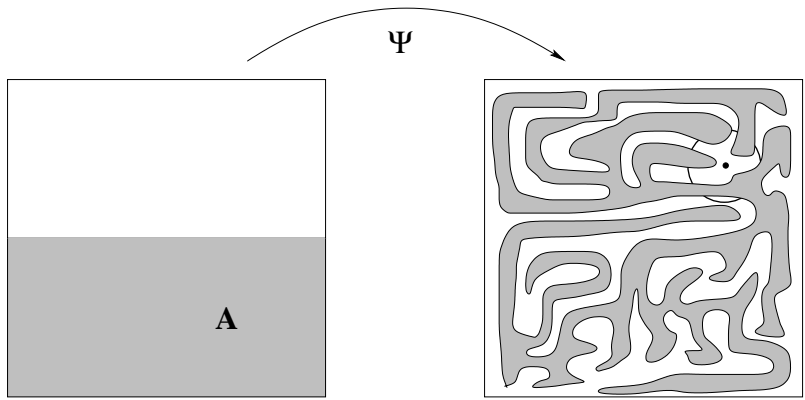

figure 2

Let $f : [0, 1] \times K \mapsto \mathbb{R}^2$ be a time dependent vector field on $K$. For any initial point $y$, denote by $t \mapsto x(t) \doteq \Psi_t(y)$ the solution of the Cauchy problem

$$\dot{x}(t) = f(t, x(t)), \qquad x(0) = y. \qquad (2.1)$$

For notational convenience, let $\widehat{\Psi} \doteq \Psi_1$ be the value of the flow map at time $t = 1$. For a fixed $0 < \kappa < 1/2$, we say that $\widehat{\Psi}$ *mixes the set $A$ up to scale $\varepsilon$* if, for every ball $B^x$ centered at a point $x \in K$ with radius $\varepsilon$, one has

$$\kappa \cdot \text{meas}(B^x) \leq \text{meas}\big(B^x \cap \widehat{\Psi}(A)\big) \leq (1 - \kappa) \cdot \text{meas}(B^x).$$

Intuitively, if a vector field has a strong mixing effect, it must be far from constant. Indeed, it should steer nearby points to quite different locations. The previous lemma thus motivates the following

**Conjecture 1.** Let $f = f(t, x)$ be a smooth vector field on $K$, and assume that the associated flow $\Psi_t$ is nearly incompressible so that, for some $\kappa' > 0$,

$$\kappa' \text{meas}(\Omega) \leq \text{meas}(\Psi_t(\Omega)) \leq \frac{1}{\kappa'} \text{meas}(\Omega), \qquad (2.2)$$

for all $\Omega \subset K$ and $t \in [0, 1]$. Then there exists a constant $C$ depending only on $\kappa, \kappa'$ such that, if $f$ mixes the set $A$ up to scale $\varepsilon$, then

$$\int_0^1 \int_K |\nabla_x f| \, dx dt \geq C \, |\log \varepsilon|. \qquad (2.3)$$

Reversing time, the above conjecture says that, in order to rearrange the set $A' \doteq \widehat{\Psi}(A)$, one needs a vector field whose cost (measured by the total variation) grows logarithmically with $\varepsilon$. A closely related conjecture can be stated as follows.

**Conjecture 2.** Consider two sets $X, Y \subset \mathbb{R}^n$ with unit measure. Assume that they are $\varepsilon$-close, in the sense that there exists a measure preserving transformation $\varphi : X \mapsto Y$ such that

$$|x - \varphi(x)| \leq \varepsilon \qquad \text{for all } x \in X.$$

Let $f = f(t, x)$ be a smooth vector field, whose flow $(t, x) \mapsto \Psi_t(x)$ is nearly incompressible and eventually separates $X$ from $Y$. More precisely we assume that (2.2) holds and that, at some time $T > 0$,

$$|\Psi_T(x) - \Psi_T(y)| \geq 1 \qquad \text{for all } x \in X, \ y \in Y.$$

Then the total variation of $f$ must be large:

$$\int_0^T \int_{\mathbb{R}^n} |\nabla_x f| \, dx dt \geq C \, |\log \varepsilon|. \qquad (2.4)$$

## Reference

[1] A. Bressan, An ill posed Cauchy problem for a hyperbolic system in two space dimensions, Preprint S.I.S.S.A., Trieste 2003.